\documentclass[12pt]{article}
\usepackage{graphicx}
\usepackage{amsmath}
\usepackage{graphics}
\usepackage{epsfig}
\usepackage{amssymb}
\usepackage{color}
 \usepackage{setspace}
\usepackage{extarrows}
\usepackage{multirow}

\usepackage[authoryear]{natbib}
\usepackage{pdflscape}
\usepackage{subfig}

\newcommand{\be}{\begin{eqnarray}}
\newcommand{\ee}{\end{eqnarray}}
\newcommand{\bea}{\begin{eqnarray*}}
\newcommand{\eea}{\end{eqnarray*}}

\newtheorem{theorem}{Theorem}[section]

\newtheorem{proposition}{Proposition}[section]
\newtheorem{corollary}{Corollary}[section]

\newtheorem{example}{Example}[section]

\newtheorem{remark}{Remark}[section]
\numberwithin{equation}{section}

\begin{document}
\title{Optimal designs for series estimation in nonparametric regression with correlated data}
\author{Kirsten Schorning, Maria Konstantinou, Holger Dette}
\author{
{\small Holger Dette, Kirsten Schorning} \\
{\small Fakult\"at f\"ur Mathematik} \\
{\small Ruhr-Universit\"at Bochum} \\
{\small 44799 Bochum, Germany} \\
\and
{\small Maria Konstantinou} \\
{\small Sustainable Energy Laboratory} \\
{\small Cyprus University of Technology} \\
{\small 3603 Limassol, Cyprus} \\
}

\maketitle
\begin{abstract}
In this paper we investigate the problem of designing experiments for series estimators in nonparametric regression models
with correlated observations. We use projection based estimators to derive an explicit solution of the best
linear oracle estimator in the continuous time model for all Markovian-type error processes. These solutions are then used to construct
estimators, which can be calculated from the available data along with their corresponding optimal design points. Our results are illustrated by means of a simulation study,
which demonstrates that the new series estimator has a better performance than the commonly used techniques
 based on the
optimal linear unbiased estimators. Moreover, we show that the performance of the estimators proposed in this paper can be further improved by
choosing the design points  appropriately.
\end{abstract}


\section{Introduction}
\label{sec1}
\def\theequation{1.\arabic{equation}}
\setcounter{equation}{0}

Nonparametric regression is a common tool of statistical inference with numerous applications [see  the monographs of
\cite{fangij1996},
  \cite{efromovich1999}, \cite{fanyao2003}, \cite{tsybakov_introduction_2009}
among many others].  The basic model is formulated in the form
\begin{equation} \label{1.1}
Y_{i} = f(X_{i}) + \varepsilon_{i}~, ~i= 1, \ldots , n,
\end{equation}
where one usually distinguishes between  random and fixed predictors  $X_{i}$. In the latter case a natural question is how to choose  $X_{1}, \ldots , X_{n}$ to obtain the most
precise estimates of the regression function $f$ and several authors have worked on this problem. For example,
 \cite{muller_optimal_1984},  \cite{biedermann_minimax_2001} and \cite{zhao_sequential_2012} derived optimal designs with respect  to different criteria  for kernel estimates,
 while  \cite{dettewiens2008} and   \cite{dette_robust_2008} considered the design problem
  for series estimation in terms of spherical harmonics and Zernike polynomials, respectively. We also refer to the work of
 \cite{efromovich_optimal_2008}, who proposed a sequential allocation scheme  in a nonparametric model of the form \eqref{1.1}
 with random predictors and heteroscedastic  errors. A common feature of the literature  in this field
 is the fact that all authors investigate the design problem in a model
 \eqref{1.1} with independent errors. However, there are many situations, where this assumption is not satisfied, in particular, when the explanatory variable represents time.
  \\
The reason for this gap in the existing literature is that the design problem  for models with correlated errors (even parametric models) is substantially harder compared to the uncorrelated case.
In contrast to the  latter case, where a very  well developed and  powerful  methodology  for  the  construction  of  optimal  designs
has been established [see,  for example, the monograph of \cite{pukelsheim2006}],
optimal designs for models with   correlated observations are only  available in rare circumstances
considering parametric models  [see, for example, \cite{pazmue2001}, \cite{natsim2003}, \cite{MulStel2004}, \cite{detleopepzhi2009}, \cite{ZhDP2010,Pazman2010}, \cite{HarmanS2010},
\cite{amosalasetal2012}, \cite{stehliketal2015}, \cite{rodriguez2017}
among others].  Some general results  on optimal designs
for linear models with correlated observations can be found in the seminal work   of
\cite{sackylv1966,sackylv1968}, while more recently  in a series of papers \cite{DetPZ2012,DetPZ2016,dette_new_2017} provided a general approach for the problem
of designing  experiments  in  linear models with correlated observations by considering the problem of optimal (unbiased linear) estimation and
optimal design simultaneously. Usually, authors use asymptotic arguments  to embed the discrete (non-convex) optimization
problem in a continuous (or approximate) one.  However, unlike the uncorrelated case, in the context of correlated observations this  approach  does  not  simplify  the  problem  substantially  and  due  to  the  lack  of  convexity
the resulting approximate optimal design problems for regression models
with correlated observations  are still extremely difficult to solve.

In this paper we  consider optimal design theory for series estimation in the nonparametric regression  model \eqref{1.1} with correlated data.
The basic notation and the general design problem are introduced in  Section \ref{sec2}.  In order to address the particular difficulties in design problems  for series estimation from  correlated data, in Section \ref{sec3}
 we consider a continuous time version of the discrete model. We first determine
optimal oracle estimators for the coefficients in a Fourier expansion of the regression function $f$. These are shrinkage
estimators  and not unbiased.  
  \\
 Section \ref{sec4} is devoted to the implementation of the results from Section \ref{sec3} for the construction of an efficient estimator with a corresponding optimal design. In particular,  we determine
an optimal approximation of the Fourier coefficients  in the continuous model  (which requires  the full trajectory of the process) by
an estimator which can be calculated from the available data $\{  Y_{t_{1}} , \ldots , Y_{t_{n}}\}$ and determine the  designs points $t_{1}, \ldots , t_{n}$
such that the approximation has minimal  mean squared error with respect to the solution in the continuous time model. The resulting estimator is a two stage estimator
shrinking the best  linear unbiased estimator when the design  points  are chosen in an optimal way.
The superiority of our approach  is demonstrated in Section \ref{sec5} by means of a  small  simulation study, while all technical details are given in Section \ref{sec6}.

\section{Optimal designs for series estimation}
\label{sec2}
\def\theequation{2.\arabic{equation}}
\setcounter{equation}{0}

Throughout this paper we consider the nonparametric regression model with a fixed design, that is,
\begin{equation}\label{disc_mod}
Y_{t_i} = f(t_i) + \varepsilon_{t_i}, \quad i=1, \ldots, n,
\end{equation}
where $f: [0,1] \to  \mathbb{R}$  is the regression function,
 $0 \leq t_1 < t_2 < \ldots < t_n\leq 1$ are $n$  distinct time points in the interval $[0, 1]$,
  $     \mathbb{E}[\varepsilon(t_j)]=0 $ and    $K(t_i,t_j) = \mathbb{E}[\varepsilon_{t_i} \varepsilon_{t_j}]$ denotes the covariance between observations at the points
   $t_i$ and $t_j$ ($ i,j=1,\dots,n $).
Let
 $$
 L^2([0, 1])= \Big\{g: [0,1]  \rightarrow \mathbb{R}: \int_0^1 g^2(t)dt < \infty  \Big\},
 $$
denote the space of square integrable (real valued)  functions with  inner product  $\langle g_{1} , g_{2} \rangle = \int_{0}^1 g_{1} (t) g_{2}(t) dt $
and norm  $\|g\|_2 = \big(\int_0^1 g^2(t) dt \big)^{1/2}$. Let $\{\varphi_j(\cdot) : j\in \mathbb{N} \}$ be an
 orthonormal basis, then any   function $f\in L^2([0, 1])$ admits a series expansion of the form
\begin{equation}\label{full_series_expan}
f(t) = \sum_{j\in \mathbb{N}} \theta_j \varphi_j(t),
\end{equation}
in $L^2{([0,1])}$
with Fourier coefficients
\begin{equation}\label{series_coef}
\theta_j = \langle f , \varphi_{j} \rangle =
\int_{0}^1 f(t) \varphi_j(t) dt \quad j \in \mathbb{N}.
\end{equation}
Moreover, the coefficients are squared summable, that is, $\sum_{j\in \mathbb{N}} \theta_j^2 < \infty$.  In order to estimate the unknown function $f$ we now follow the idea of projection estimators [see \cite{tsybakov_introduction_2009}, pp.47]
and estimate the  truncated series  ${f}^{(J)}(t) = \sum_{j=1}^{J} \theta_j \varphi_j(t) $ by
\begin{equation}\label{h3}
\hat{f}^{(J)}(t) = \sum_{j=1}^{J} \hat\theta_j \varphi_j(t),
\end{equation}
where $ \hat\theta_j$ is an appropriate estimator for the Fourier coefficient $ \theta_j $ $(j  = 1 ,  \ldots , J$). For example, if $\max^n_{i=2} (t_i - t_{i-1}) \to 0$, as $n \to \infty$, an asymptotically unbiased estimator  of  $\theta_j$
is given by
\begin{equation} \label{riem}
\sum_{i=2}^{n} (t_{i} -t_{i-1}) \varphi_j ({t_{i-1}})  Y_{t_{i-1}}.
\end{equation}
More general estimators will be specified later on.
At this point it is only important to  note that the performance  of any reasonable estimator will depend on the design points $t_{1},\ldots , t_{n}$.
We are interested in  choosing these design points such that
the mean integrated squared error
$$
\mathbb{E} \Big [ \int_{0}^{1}   \big ( \hat{f}^{(J)}(t)  - f(t) \big  )^2   \Big]  =   \sum_{j=1}^{J} \mathbb{E} \big [ ( \hat \theta_{j} - \theta_{j} )^{2} \big ]
+   \sum_{j=J+1}^{\infty} \theta_{j}^{2},
$$
is minimal.  We also note that  any solution of this discrete optimization problem depends on the unknown regression function $f$, the truncation point $J$ used in  \eqref{h3} and
on the covariance kernel $K$, which is assumed to be known throughout this paper.  On the other hand, the term $ \sum_{j=J+1}^{\infty} \theta_{j}^{2}$ does not depend on the design points
which can therefore, be determined by minimizing $  \sum_{j=1}^{J} \mathbb{E} \big [ ( \hat \theta_{j} - \theta_{j} )^{2} \big ]$ with respect to the choice of $t_1, \ldots, t_n$.
For example, if $\hat \theta_{j} =  \sum_{i=1}^{\ell_{j}}\alpha_{ji} Y_{t_{i}}$ is a linear estimator of $\theta_{j}$ ($j=1,\ldots , J$) we
have that
\begin{eqnarray} \label{discoptdes}
 \sum_{j=1}^{J} \mathbb{E} \big [ ( \hat \theta_{j} - \theta_{j} )^{2} \big ] &=& \sum_{j=1}^{J} \Big (\sum_{i=1}^{\ell_{j}} \alpha_{ji} f(t_{i}) - \theta_{j }\Big ) ^2
 +  \sum_{j=1}^{J} \sum_{i_{1},i_{2}=1}^{\ell_{j}}  \alpha_{ji_{1}} \alpha_{ji_{2}} K(t_{i_{1}},t_{i_{2}}),
\end{eqnarray}
which has to be minimized with respect to the choice of the time points $t_{1}, \ldots , t_{n}$.

\section{Optimal estimation in the continuous time model}
\label{sec3}
\def\theequation{3.\arabic{equation}}
\setcounter{equation}{0}
The discrete optimization  problem \eqref{discoptdes} stated in the previous section is extremely difficult to be solved.
In this section in order to derive efficient designs, we investigate a simpler problem and
 consider the continuous time  nonparametric regression model of the form
\begin{equation}\label{cont_mod}
Y_t = f(t) + \varepsilon_t \, , \quad t \in [0, 1],
\end{equation}
where $f $ is an unknown square  integrable   function and the error process $\varepsilon= \{\varepsilon_t : t \in [0, 1]\}$ is a
centered Gaussian process with covariance kernel $K(s,t)=\mathbb{E} [\varepsilon_s \varepsilon_t]$.
As we assume that the full trajectory of the  process is  available, there is in fact no optimal design problem but only the issue of optimal estimation of the regression function  $f$.
The optimal design question will appear later, when we return to the discrete model \eqref{disc_mod}.
The main result of this section provides an oracle  solution of the optimal estimation problem. In particular, the optimal estimator depends on the unknown function $f$ in  model \eqref{cont_mod} and is therefore,
not implementable (even if the full trajectory of the process $\{ Y_t : t \in [0,1]  \} $ is   available). However, our solution serves as benchmark
and actually provides a clear hint how good estimators and corresponding optimal designs can be constructed. This will be formulated precisely in  Section \ref{sec4}.

Model \eqref{cont_mod}  is often written in terms of a stochastic
differential equation (provided that the regression function $f$ is differentiable with derivative $\dot{f}$), that is
\begin{equation}\label{cont_modif}
dY_t = \dot{f}(t) dt + d  \varepsilon_t \, , \quad t \in [0, 1]~,
\end{equation}
If $\varepsilon=\{\varepsilon_t : t \in [0, 1]\}$ is a Brownian motion, the model \eqref{cont_modif}  is called Gaussian white noise model
and has found much attention in the statistical literature [see, for example, \cite{MR620321} or  \cite{tsybakov_introduction_2009}  among many others].
In particular, the model
is asymptotically equivalent  to   the nonparametric regression model $Z_{i} = \dot{f}(i/n)  +   \eta_{i }$  ($i=1 , \ldots ,n$), where
$\eta_{1 }, \ldots , \eta_{n}$ are independent standard normally distributed random variables
  [see  \cite{brown1996}].
Note that the focus in the aforementioned publications is on the optimal  estimation of the function $ \dot{f}$, whereas in this section we are interested in the estimation the function $f$
in model  \eqref{cont_mod}.
Nevertheless, under additional assumptions we can investigate the properties of the derivative of the oracle estimator developed in what follows  
and a brief discussion of these relations is given   in Example \ref{exam_cont_modif}. \\
Another important difference between model \eqref{cont_mod}  and the Gaussian white noise model commonly discussed in the literature of mathematical statistics lies in the fact that we consider a general
 error process  $\{\varepsilon_t : t \in [0, 1]\}$. In particular, we concentrate on  Markovian Gaussian error processes
with a covariance kernel of the form
\begin{equation}\label{triangular-kernel}
\mathbb{E}[\varepsilon_s\varepsilon_t]= K(s, t) = u(s) v(t) \quad \mbox{for } s \leq t,
\end{equation}
where  $u(\cdot)$ and $v(\cdot)$ are some (known) functions defined on the interval $[0, 1]$, such that $v(t) \neq 0$ for $t\in [0, 1]$. Kernels of this form generalize
 the Brownian motion, which is obtained for $u(t) =t$, $v(t)=1$, and are called triangular kernels in the literature. The property \eqref{triangular-kernel} essentially characterizes a Gaussian process to be Markovian [see \cite{doob_heuristic_1949} or \cite{mehr1965certain} for more details].
We assume that the process $\{\varepsilon_t : t\in [0,1]\}$ is non-degenerate on the open interval $(0,1)$, which implies that the function
\begin{equation}\label{qfct}
q(t) = \frac{u(t)}{v(t)},
\end{equation}
is positive on the interval $(0,1)$ and strictly increasing and continuous on $[0,1]$.

Regarding the estimation of  the unknown function $f$, we propose to  estimate the coefficients $\theta_j$ in the  projection estimator
 \eqref{full_series_expan} using statistics of the form  [see \cite{grenander1950}]
  \begin{equation}\label{linear_estimate}
\hat\theta_j= \int_{0}^{1} Y_{t} \xi_j(dt), ~ j \in \mathbb{N} ,
\end{equation}
where $\xi_j$ is a signed measure on the interval $[0,1]$ such that
\begin{equation}\label{h1}
\sum_{j=1}^{\infty}  \big \{  (\mathbb{E} [\hat \theta_{j }] )^2  + \mbox{Var} (\hat \theta_{j})  \big \} =
\sum_{j=1}^{\infty} \Big(\int_0^1 f(t) d\xi_j(t) \Big)^2 + \sum_{j=1}^{\infty} \int^1_0 \int^1_0 K(s, t) d\xi_j(s) d\xi_j(t) < \infty.
\end{equation}
Obviously,  this condition implies for the sequence of estimators $(\hat\theta_j)_{j\in \mathbb{N}}$ that
$\sum_{j=1}^{\infty} \mathbb{E}[\hat\theta^2_j] < \infty , $
and thus we can define the random variable
\begin{equation}\label{h2}
\hat f (t) = \sum_{j=1}^{\infty} \hat\theta_j \varphi_j(t).
\end{equation}
In particular,  if $
\hat{f}^{(J)}(t) = \sum_{j=1}^{J} \hat\theta_j \varphi_j(t)
$
is the truncated series from  \eqref{h2},  we have that
\begin{eqnarray*}
\lim_{J \rightarrow \infty} \mathbb{E}\Big[\int_0^1 \big(\hat f^{(J)}(t) - f(t) \big)^2 dt \Big]
=  \lim_{J\rightarrow \infty} \sum_{j=1}^{J} \mathbb{E}[(\hat\theta_j - \theta_j)^2]
= \sum_{j=1}^{\infty} \mathbb{E}[(\hat\theta_j - \theta_j)^2] < \infty,
\end{eqnarray*}
and  the mean integrated squared error of the estimator $\hat f$ in \eqref{h2} is given by
\begin{eqnarray}\label{miseinfinity}
\mbox{MISE}(\hat f)& :=& \mathbb{E} \Big [\int_{0}^1 (\hat f (t) - f(t))^2 dt \Big ]  = \sum_{j=1}^\infty \mathbb{E}[(\hat\theta_j - \theta_j)^2].
\end{eqnarray}

\medskip

We conclude that  the optimal linear oracle estimator $\hat f$ of the function $f$  minimizing
 \eqref{miseinfinity} can be determined  minimizing
the individual mean squared errors  $\mathbb{E}[(\hat\theta_j - \theta_j)^2]$  separately.  Due to the definition of linear estimators in \eqref{linear_estimate},
this problem  corresponds to the determination of a signed measure $\xi^*_j$ on the interval $[0, 1]$, which minimizes the functional
\begin{equation}\label{mse_criterion}
\begin{split}
\Psi_j(\xi_j) :&= \mathbb{E}\Big[\Big( \int_{0}^{1} Y_{t} \xi_j(dt)- \theta_j\Big)^2\Big]\\
&= \int_0^1 \int_0^1 \big [  f(s)f(t) + K(s, t)  \big ] \xi_j(ds) \xi_j(dt) - 2\theta_j \int_0^1 f(s) \xi_j(ds) + \theta^2_j \\
&= \int_0^1 \int_0^1  K(s, t) \xi_j(ds) \xi_j(dt) + \Big ( \int_0^1 f(s) \xi_j(ds) - \theta_j\Big)^2 \, .
\end{split}
\end{equation}

\begin{remark} \label{remgren}
  {\rm~~
  \begin{itemize}
  \item[(1)]
 Note that - in contrast to most of the literature - we do not assume that $\hat\theta_j$ is  an unbiased  estimator  of the Fourier coefficient $\theta_j$ $ (j \in \mathbb{N})$.
 A prominent unbiased estimator for $\theta_j$ is given by
 \begin{equation}\label{h4}
 \tilde{\theta}_j = \int_0^1 Y_t \varphi_j(t) dt \quad (j\in \mathbb{N}) \, ,
 \end{equation}
 and for general unbiased estimates of the form \eqref{linear_estimate} the condition \eqref{h1} reduces to
 \begin{equation} \label{h5}
 \sum_{j=1}^{\infty} \int^1_0 \int^1_0 K(s, t) d\xi_j(s) d\xi_j(t) < \infty \, .
 \end{equation}
Moreover, if the kernel $K$ is continuous on $[0, 1]\times [0, 1]$ and if $\varphi_1, \varphi_2, \ldots$ are the eigenfunctions of the integral operator associated with
 the covariance kernel $K$ with corresponding eigenvalues $\lambda_1, \lambda_2, \ldots$, then  condition  \eqref{h1} further reduces to
 \begin{eqnarray*}
 \sum_{j=1}^{\infty} \int_0^1 \int_0^1 K(s, t) \varphi_j(s) \varphi_j(t) ds dt = \sum_{j=1}^\infty \lambda_i \int_0^1 \varphi_j(t) \varphi_j(t) dt = \sum_{j=1}^\infty \lambda_j < \infty.
 \end{eqnarray*}
  \item[(2)]
  Under the additional assumption that the estimator  \eqref{linear_estimate} is unbiased for  $\theta_{j}$, the second term
in \eqref{mse_criterion} vanishes and the resulting optimization problem corresponds to the problem of finding
the best linear estimator in the location scale model $Y_{t} =  \theta_{j} + \varepsilon_{t}$, which has been first studied
  in a seminal paper of \cite{grenander1950}. This author showed that under the additional constraint $\int_{0}^{1} d \xi_{j}  (dt) =1 $
  the optimal solution $\xi_{j}^{*}$ minimizing $\int_0^1 \int_0^1  K(s, t) \xi_j(ds) \xi_j(dt) $  can be characterized by the
  property that  the function  $t \to \int_0^1  K(s, t) \xi_j^{*}(ds)  $ is constant on the interval $[0,1]$.
  \end{itemize}
  }
\end{remark}

 The following theorem provides a complete solution of the optimization problem \eqref{mse_criterion}
  and is proven in the appendix.
 For a precise statement of the result we denote by $\delta_x$ the Dirac measure at the point $x$ and distinguish the following cases
 for the triangular kernel  \eqref{triangular-kernel}.
 \begin{itemize}
 \item[(A)] $u(0) \neq 0$.
 \item[(B)]$u(0) = 0, f(0) = 0$.
 \item[(C)] $u(0) = 0 , f(0) \neq 0 $.
 \end{itemize}

\begin{theorem}\label{theo_cont_opt_est}
Consider the functional   $\Psi_j$ in \eqref{mse_criterion}  with a twice differentiable regression function $f$ and 
 a triangular covariance kernel of the form \eqref{triangular-kernel}, where the functions $u$ and $v$ are  also twice differentiable.
For any $j \in  \mathbb{N} $  the signed measure $\xi^*_j(dt)$
 minimizing the functional $\Psi_j$  in the class of all signed measures on the interval
 $[0,1]$  is given by
\begin{equation}\label{opt_cont_measure}
\xi^*_j(dt) = \frac{\theta_j}{1+ c} \left(P_0 \delta_0(dt) + P_1\delta_1(dt) + p(t) dt \right),
\end{equation}
where $\theta_j$ is the $j$-th Fourier coefficient in the Fourier expansion \eqref{full_series_expan}.
The values for $c$, $P_0$, $P_1$ and the function $p(\cdot)$ do not depend on the index $j$ and take different values corresponding to the properties of the functions $u(\cdot)$ and $f(\cdot)$. In particular, we have the following cases
\begin{enumerate}
\item[(A)] If $u(0)\neq 0$, the quantities $c$, $P_0$, $P_1$ and $p$ are given by
\begin{eqnarray}
c&=& \int_{0}^{1} \Big\{\frac{d}{dt} \left[ \frac{f(t)}{v(t)} \right]\Big\}^2 \Big(\frac{d}{dt}{q(t) }\Big)^{-1} dt + \frac{f^2(0)}{v^2(0)}{(q(0))^{-1}}, \label{c}\\
P_0 &=& - \frac{1}{v(0)} \frac{d}{dt}\left[\frac{f(t)}{u(t)}\right]\Big |_{t=0}  \Big(\frac{d}{dt}{q(t) }\Big |_{t=0} \Big)^{-1}q(0),  \label{P0}\\
P_1&=& \frac{1}{u(1)} \frac{d}{dt}\left[\frac{f(t)}{v(t)}\right]\Big |_{t=1}  \Big(\frac{d}{dt}{q(t) }\Big |_{t=1} \Big)^{-1}q(1),  \label{P1} \\
p(t) &=&- \frac{1}{v(t)}\frac{d}{dt}\Big\{\frac{d}{dt} \left[ \frac{f(t)}{v(t)} \right] \Big(\frac{d}{dt}{q(t) }\Big)^{-1}\Big\}  \label{pt},
\end{eqnarray}
where the function $q$ is defined in \eqref{qfct}.
\item[(B)] If $u(0) = 0$ and $f(0)=0$, the quantities $c$ and $P_0$ are given by
\begin{eqnarray}
c&=& \int_{0}^{1} \Big\{\frac{d}{dt} \left[ \frac{f(t)}{v(t)}\right]\Big\}^2 \Big(\frac{d}{dt}{q(t) }\Big)^{-1} dt , \label{c2}\\
P_0 &=& 0  \label{P02}
\end{eqnarray}
and $P_1$ and $p$ are given by \eqref{P1} and \eqref{pt}, respectively.
\item[(C)]If $u(0) = 0$ and $f(0) \neq 0$, the quantities $c$, $p(t)$ and $P_1$ are equal to zero, whereas $P_0$ is given by
\begin{equation}
P_0 = \frac{1}{f(0)}.
\end{equation}
\end{enumerate}
\end{theorem}

\begin{corollary}\label{corr_cont_opt_est}
Consider the regression model \eqref{cont_mod} with a twice differentiable regression function $f$ and a non-degenerate centered Gaussian error process $\{\varepsilon_t: t \in [0, 1]\}$ with a triangular covariance kernel of the form \eqref{triangular-kernel}, where the functions $u$ and $v$ are twice differentiable.
The best linear oracle estimator minimizing the mean integrated squared error in \eqref{miseinfinity} in the class of all linear estimators of the form \eqref{h2}
satisfying \eqref{h1} is defined by  
$$f^* (t)  = \sum_{j=1}^\infty  \hat\theta_j^* \varphi_j (t),$$
 where the coefficients  $\hat\theta_j^*$ are given by
\begin{equation} \label{optest}
\hat\theta_j^*= \int_0^1 Y_t \xi^*_j(dt)~,~~
j \in \mathbb{N} ,
\end{equation}
and  the signed measure $\xi^*_j(dt)$ is
defined in Theorem \ref{theo_cont_opt_est}.
Moreover, the corresponding mean integrated squared error is given by
$$
\mbox{\rm MISE}(\hat f^*) =  \frac{1}{1+c} \sum_{j=1}^\infty \theta_j^2 =    \frac{1}{1+c}  \int_0^1 f^2(t) dt \ ,
$$
where $c$ is defined in \eqref{c}.
\end{corollary}
{\rm  Note that Theorem \ref{theo_cont_opt_est} is a theoretical result as it requires knowledge of the unknown regression function $f$.
Nevertheless, we will use it extensively in the following section to construct good estimators and corresponding optimal designs for series estimation in model \eqref{disc_mod}.}

\begin{remark} {\rm  ~~
\begin{itemize}
\item[(1) ]
In model \eqref{disc_mod} with covariance kernel \eqref{triangular-kernel} and $u(0)=0$, the observation $Y_0$ at $t=0$ does not contain any error.
Therefore, the value of $f(0)$ is known so that  it can be checked whether case (B) or (C) of Theorem \ref{theo_cont_opt_est} holds.
\item[(2)] The estimator given in Theorem \ref{theo_cont_opt_est} depends on the orthonormal system of the series expansion via the parameter $\theta_j$.
\item[(3)]  Using integration by parts the resulting estimator $\hat\theta^*_j$ in Theorem \ref{theo_cont_opt_est} can be represented as stochastic integral. For example, in case (A)
 (where $u(0) \neq 0$) the estimator can be represented as
\begin{equation}\label{reformulated_est}
(A) ~~~~~~~~~~~~~~~
\hat\theta^*_j= \frac{\theta_j}{1+c}\Big\{ \int_{0}^1 {\frac{d}{dt}\Big[\frac{f(t)}{v(t)}\Big]}  \left(\frac{d}{dt}{q(t) }\right)^{-1} d\Big(\frac{Y_t}{v(t)} \Big) + \frac{f(0)}{u(0)} \frac{Y_0}{v(0)}\Big\},
\end{equation}
where the constant $c$ is defined in \eqref{c}.  Similarly in case (B) (where $u(0)=0$ and $f(0)=0$), the estimator can be represented by
\begin{equation}\label{reformulated_estu0}
(B) ~~~~~~~~~~~~~~~  \hat\theta^*_j= \frac{\theta_j}{1+c}\Big\{ \int_{0}^1 {\frac{d}{dt}\Big[\frac{f(t)}{v(t)}\Big]}  \left(\frac{d}{dt}{q(t) }\right)^{-1} d\Big(\frac{Y_t}{v(t)} \Big)\Big\}.
~~~~~~~~~~~~~~~
\end{equation}
Finally, in case (C) (where $u(0)=0$ and $f(0)\neq 0$), the estimator directly reduces to
\begin{equation}\label{reformulated_est3}
(C) ~~~~~~~~~~~~~~~   \hat\theta^*_j= \theta_j. ~~~~~~~~~~~~~~~ ~~~~~~~~~~~~~~~ ~~~~~~~~~~~~~~~ ~~~~~~~~~~~~~~~  ~~~~~~~
\end{equation}
In the latter case the estimator in \eqref{reformulated_est3} is not random, but fixed to the true - but unknown - parameter  $\theta_j$.
\end{itemize}
}
\end{remark}

\begin{example}\label{exam_cont_modif}
{\rm
A very popular orthonormal basis of $L^2([0,1])$ is given by the trigonometric functions
\begin{equation}\label{cosine_basis}
\varphi_j(t) = \begin{cases}
1 \quad, & j=1 \\
\sqrt{2} \cos(2\pi k t) \quad, &j=2k \\
\sqrt{2} \sin(2\pi k t) \quad, & j= 2k+1
\end{cases}, \quad j=1, 2, \ldots \ .
\end{equation}
Under the assumptions of Theorem \ref{theo_cont_opt_est} we assume  that $f$ and its derivative $\dot{f}$ can be represented as a trigonometric series, that is,
\begin{eqnarray}
\label{trig1}
f(t)   &=&  \theta_1  + \sum_{k=1}^{\infty} \sqrt{2}\cos(2\pi k t) \theta_{2k} + \sum_{k=1}^{\infty} \sqrt{2}\sin(2\pi k t) \theta_{2k+1} ~, \\
\dot{f} (t ) &=&
\bar{\theta}_1  + \sum_{k=1}^{\infty} \sqrt{2}\cos(2\pi k t) \bar{\theta}_{2k} + \sum_{k=1}^{\infty} \sqrt{2}\sin(2\pi k t)  \bar{\theta}_{2k+1}~.
\label{trig2}
\end{eqnarray}
Note that  (under suitable assumptions) the Fourier coefficients  in \eqref{trig1}
and  \eqref{trig2}  are related by the equations
\begin{equation} \label{relation_coeffs}
\bar{\theta}_1 = 0 , \quad \bar{\theta}_{2k} =  (2\pi k) \theta_{2k+1},   \quad \bar{\theta}_{2k+1}=  -(2\pi k) \theta_{2k}.
\end{equation}
If the error process $\{\varepsilon_t: t \in [0, 1]\} $  in model \eqref{disc_mod} is given by a  Brownian motion, we have  $u(t) =t$, $v(t) =1$ in the definition of the triangular kernel  \eqref{triangular-kernel} and thus $q(t) = t$.
A straightforward application of  Corollary \ref{corr_cont_opt_est} (case (B)) yields for the optimal oracle estimator of the function $f$
\begin{equation}\label{sincos_est}
f^*(t) = \hat\theta^*_1 + \sum_{k=1}^{\infty} \sqrt{2}\cos(2\pi k t) \hat\theta^*_{2k} + \sum_{k=1}^{\infty} \sqrt{2}\sin(2\pi k t) \hat\theta^*_{2k+1}\ ,
\end{equation}
 where the estimated Fourier coefficients are given by
\begin{equation}\label{best_fourier_est}
\hat\theta^*_j = \frac{\theta_j}{1+c} \int_0^1 \dot{f}(t) dY_t  ~,~j \in \mathbb{N}~,
\end{equation}
(note that $f(0) = f(1)=0$).
We thus also   obtain an  estimator of  the function $\dot{f}$ in model
\eqref{cont_modif} by taking the derivative of $f^*$ given n \eqref{sincos_est}, that is,
\begin{equation}\label{fhat_dot}
\dot{f^*}(t) = - \sum_{k=1}^{\infty} (2\pi k) \sqrt{2}\sin(2\pi k t) \hat\theta^*_{2k} + \sum_{k=1}^{\infty} (2\pi k) \sqrt{2}\cos(2\pi k t) \hat\theta^*_{2k+1} \ .
\end{equation}
Using the relation \eqref{relation_coeffs},  the estimator in \eqref{best_fourier_est} can be rewritten as
\begin{equation*}
\hat\theta^*_j= \begin{cases}
-\frac{\bar{\theta}_{2k+1}}{2\pi k } \frac{1} {1+c} \int_0^1 \dot{f}(t) dY_t \  ,  \quad & j= 2k \\
\frac{\bar{\theta}_{2k}}{2\pi k } \frac{1} {1+c} \int_0^1 \dot{f}(t) dY_t \ , \quad & j= 2k +1 ,
 \end{cases}
\end{equation*}
and the mean integrated squared error of the estimator $\dot{f^*}$ in \eqref{fhat_dot} is given by
\begin{eqnarray} \label{mise_fhatdot}
\mathbb{E}\Big[\int_0^1 \big(\dot{f^*}(t) - \dot{f}(t)  \big)^2 dt \Big]&=& \sum_{j=2}^{\infty}\frac{\bar{\theta}^2_{j}}{(1+c)^2} \mathbb{E}\Big[\big(1+c - \int_{0}^{1} \dot{f}(t) dY_t\big)^2\Big] \\ 
&=& \notag \sum_{j=2}^{\infty} \frac{\bar{\theta}^2_{j}}{1+c}= \frac{\sum_{j=2}^{\infty} \bar{\theta}^2_{j}}{1+\sum_{j=2}^{\infty} \bar{\theta}^2_{j} },
\end{eqnarray}
where we have used  the representation  $c= \int_{0}^{1}\big( \dot{f}(t) \big)^2 dt = \sum_{j=1}^{\infty} \bar{\theta}^2_{j}= \sum_{j=2}^{\infty} \bar{\theta}^2_{j}$
in the last equality. \\
It might be of interest to compare this  estimator with the linear oracle estimator
\begin{equation}\label{oracle}
\dot{\tilde{f}}(t) = \sum_{j\in \mathbb{N}} \tilde{\theta}_j \varphi_j(t),
\end{equation}
proposed in  \cite{tsybakov_introduction_2009}[p. 67], where
$$\tilde{\theta}_j = \frac{\bar{\theta}^2_j}{1 + \bar{\theta}^2_j}\int_0^1 \varphi_j(t) dY_t \,  ,$$
is used as the estimator of the Fourier coefficient $\bar{\theta}_j$ ($j=1, 2, \ldots ).$  This estimator is a shrinkage version of the unbiased estimator in \eqref{h4}
and the mean integrated  squared error of $\dot{\tilde{f}}$ is given by
\begin{eqnarray}\label{mise_ftildedot}
\mathbb{E}\Big[\int_0^1 \big(\dot{\tilde{f}}(t) - \dot{f}(t)  \big)^2 dt \Big] &=& \sum_{j=1}^{\infty} \frac{\bar{\theta}^2_j}{1+ \bar{\theta}_j^2} .
\end{eqnarray}
Comparing \eqref{mise_fhatdot} and \eqref{mise_ftildedot}, we observe that the oracle estimator $\dot{f^*}$, which is constructed by an application of Corollary \ref{corr_cont_opt_est}, has a smaller mean integrated squared error than the   estimator $\dot{\tilde{f}}$ defined in \eqref{oracle}.
}
\end{example}

 \section{Efficient  series estimation from  correlated data}
 \label{sec4}
 \def\theequation{4.\arabic{equation}}
\setcounter{equation}{0}

 In this section we apply the results from the continuous time model to construct optimal designs for series estimation of  the function $f$
 in model \eqref{disc_mod}.  In  this  transition from the continuous to the discrete model we are faced with several challenges. First,
  the signed measure defining the optimal oracle estimator $\hat\theta_j^*$  depends on the unknown function $f$
  through its Fourier coefficients  and through the constant $c$, and  the function $f$ also appears in the stochastic integrals
   in \eqref{reformulated_est} and \eqref{reformulated_estu0}.  Secondly, we need  to address the problem that even
   with preliminary knowledge of the
   function $f$, the stochastic integrals can not be computed since as the continuous time process $\{ Y_t : t \in [0,1] \} $ is not observable.
   In order to overcome these difficulties and construct an implementable  estimator, which does not require preliminary knowledge of $f$, we proceed  to several steps, which are explained in detail below. Roughly speaking,  these steps consist of a two stage estimation procedure, a truncation and an
   appropriate approximation of the stochastic integrals by sums, which can be calculated from the available data.
   In the latter step of this procedure we also determine the optimal design points.  \\
 Throughout this section  we will restrict ourselves to the cases (A) and (B) of Theorem \ref{theo_cont_opt_est}.
For the case  (C) we simply propose to replace the parameter value \eqref{reformulated_est3} by the best linear unbiased estimator derived
in \cite{dette_new_2017}.
\subsection{Truncation in the continuous time model}
 \label{sec41}

In model \eqref{disc_mod}  with $n$ observations, only a finite number, say $J$, of Fourier coefficients in the series expansion \eqref{full_series_expan} can be estimated. For this reason, we consider
for fixed $J\in \mathbb{N}$
the best $L^2$-approximation
 \begin{equation}\label{f_truncated}
f^{(J)}(t)  =  \sum_{j=1}^J \theta_j \varphi_j(t) = \Phi^{(J), T}(t) \theta^{(J)},
\end{equation}
 of the function $f$ by functions from the $\mbox{span}\{\varphi_1, \ldots, \varphi_J\}$, space 
where  the vectors    $\theta^{(J) }$  and  $\Phi^{(J)} $  are defined  by    $\theta^{(J)}=(\theta_1, \ldots, \theta_J)^T$ and $\Phi^{(J)}(t) = (\varphi_1(t), \ldots, \varphi_J(t))^T$,
respectively. We now replace the function $f$ by the function $f^{(J)}$ in the  estimators  $\hat \theta_1^*, \ldots, \hat  \theta_J^*$  defined in \eqref{reformulated_est} and \eqref{reformulated_estu0} for cases (A) and (B) respectively. In case (A) this gives the vector
 \begin{equation}\label{truncated_estuA}
\hat\theta^{(J),*}= \frac{1}{1+c^{(J)}}\theta^{(J)}(\theta^{(J)})^{T}\Big\{ \int_{0}^1 \frac{d}{dt}\Big[\frac{\Phi^{(J)}(t)}{v(t)}\Big]\left(\frac{d}{dt}{q(t) }\right)^{-1} d\Big(\frac{Y_t}{v(t)} \Big)+ \frac{\Phi^{(J)}(0)}{u(0)} \frac{Y_0}{v(0)}\Big\}  \quad ,
\end{equation}
where
\begin{eqnarray}\label{cJ2}
c^{(J)} &=&   (\theta^{(J)})^T C^{(J)} \theta^{(J)}, \label{cJ2}
\end{eqnarray}
and the $J \times J$ matrix $ C^{(J)}$  is defined by
 \begin{equation}\label{CJmat}
 C^{(J)} =  \int_0^1\frac{d}{dt}\Big[\frac{\Phi^{(J)}(t)}{v(t)}\Big] \Big(\frac{d}{dt}\Big[\frac{\Phi^{(J)}(t)}{v(t)}\Big]\Big)^T \Big(\frac{d}{dt}{q(t) }\Big)^{-1}  dt+ \frac{\Phi^{(J)}(0) (\Phi^{(J)}(0))^T}{u(0) v(0)} ~. \end{equation}
 Similarly, in case (B) we obtain
  \begin{equation}\label{truncated_estuB}
\hat\theta^{(J),*}= \frac{1}{1+m^{(J)}}\theta^{(J)}(\theta^{(J)})^{T}\Big\{ \int_{0}^1 \frac{d}{dt}\Big[\frac{\Phi^{(J)}(t)}{v(t)}\Big]\left(\frac{d}{dt}{q(t) }\right)^{-1} d\Big(\frac{Y_t}{v(t)} \Big)\Big\} ~,
\end{equation}
where
 \begin{eqnarray}\label{cJu0}
m^{(J)}
&=&   (\theta^{(J)})^T M^{(J)} \theta^{(J)},\label{cJ2u0}
\end{eqnarray}
and the  $J \times J$ matrix $M^{(J)}$ is given by
\begin{equation}\label{MJmat}
M^{(J)} =  \int_0^1\frac{d}{dt}\Big[\frac{\Phi^{(J)}(t)}{v(t)}\Big] \Big(\frac{d}{dt}\Big[\frac{\Phi^{(J)}(t)}{v(t)}\Big]\Big)^T \Big(\frac{d}{dt}{q(t) }\Big)^{-1}  dt �.
 \end{equation}

The resulting estimators \eqref{truncated_estuA} and \eqref{truncated_estuB} still depend on the first $J$ unknown Fourier coefficients $\theta_1, \ldots, \theta_J$ and also depend on the full trajectory of the process $\{Y_t : t \in [0, 1]\}$.
This dependence will be removed in the following sections.

 \subsection{Discrete approximation of stochastic integrals}  \label{sec42}

 In concrete applications   the integrals in  \eqref{truncated_estuA} and   \eqref{truncated_estuB}  cannot
 be evaluated and have to be approximated from the given data. For this purpose
 we assume that $n$ observations $Y_{t_1}, \ldots, Y_{t_n}$ from model \eqref{disc_mod}
at $n$ distinct time points $0=t_1< t_2 < \ldots < t_{n-1} < t_n = 1$ are available  and we  consider the estimators
 \begin{eqnarray}  \label{hol10}
 \hat\theta^{(J), n}
 &=& \frac{1}{1+c^{(J)}}\theta^{(J)}(\theta^{(J)})^{T} \Big\{ \sum_{i=2}^{n}\mu_i (\frac{Y_{t_i}}{v(t_i)} - \frac{Y_{t_{i-1}}}{v(t_{i-1})} ) + \frac{\Phi^{(J)}(0)}{u(0)} \frac{Y_0}{v(0)} \Big\}, \label{lin_est_disc}
 \\
 \label{hol11}
 \hat\theta^{(J), n}
 &=& \frac{1}{1+m^{(J)}}\theta^{(J)}(\theta^{(J)})^{T} \Big\{ \sum_{i=2}^{n}\mu_i (\frac{Y_{t_i}}{v(t_i)} - \frac{Y_{t_{i-1}}}{v(t_{i-1})} ) \Big\},
 \end{eqnarray}
 as approximations of  the quantities  in \eqref{truncated_estuA}  and  \eqref{truncated_estuB}, respectively.
Note that $\hat\theta^{(J), *}$ depends on the full trajectory $\{ Y_t : t \in [0,1] \}$, while $\hat\theta^{(J), n}$ is an approximation based on the sample $\{ Y_{t_{i} }: i=1,\ldots,n \}$.
 In \eqref{lin_est_disc} and \eqref{hol11} $\mu_2, \dots , \mu_n$ denote $J$-dimensional weights   which depend on the  time points $0=t_1 <t_2 <\ldots < t_{n-1}< t_n=1$ and will be chosen in an optimal way.
 In particular we propose to determine  the weights $\mu_2, \ldots, \mu_n$   such that the expected $L^2$-distance
\begin{equation}\label{l2_festimates}
\mathbb{E}\big[\|\hat\theta^{(J), *} - \hat\theta^{(J), n} \|^2\big]
\end{equation}
 between         $\hat\theta^{(J), *}$ and its discrete analogue  $\hat\theta^{(J), n}$  is minimized, where $\|\cdot\|$ denotes the Euclidean norm in $\mathbb{R}^{J}$.
 
The following result provides an alternative expression of the expectation of this distance.
 \begin{proposition}\label{lem_reformulated_mse}
 Assume that the conditions of Theorem \ref{theo_cont_opt_est} are satisfied. The Euclidean distance between
 the estimators  $ \hat\theta^{(J), *}$ and   $\hat\theta^{(J), n}$
 can be  represented as
\begin{equation}\label{reformulated_mse}
\mathbb{E}\big[\|\hat\theta^{(J), *} - \hat\theta^{(J), n} \|^2\big] =  k^{(J)} \big \{  V (\mu_2, \ldots ,\mu_n)  + B (\mu_2, \ldots ,\mu_n)  \big \}~,
\end{equation}
where the quantities $V$ and $B$ are defined by
\begin{eqnarray} \label{hol6}
V (\mu_2, \ldots ,\mu_n)
&=&   \mbox{\rm{tr}}\Big \{\sum_{i=2}^{n}\int_{t_{i-1}}^{t_i} \Big({\frac{d}{dt}\Big[\frac{\Phi^{(J)}(t)}{v(t)}\Big]} ~  \\
&&
 \nonumber
 \times
 \Big(\frac{d}{dt}{q(t) }\Big)^{-1} - \mu_i\Big) \Big({\frac{d}{dt}\Big[\frac{\Phi^{(J)}(t)}{v(t)}\Big]}   \Big(\frac{d}{dt}{q(t) }\Big)^{-1} - \mu_i\Big)^T \Big(\frac{d}{dt}{q(t) }\Big)dt  \Big\},
  \\
 B (\mu_2, \ldots ,\mu_n) &=& \nonumber
  \mbox{\rm{tr}} \Big \{
  \sum_{i=2}^n \int_{t_{i-1}}^{t_i} \Big({\frac{d}{dt}\Big[\frac{\Phi^{(J)}(t)}{v(t)}\Big]}  \left(\frac{d}{dt}{q(t) }\right)^{-1} - \mu_i\Big)\Big(\frac{d}{dt}\Big[\frac{{f(t)}}{v(t)}\Big]\Big)dt \\\
\label{hol7}
&&
 \times \Big(\sum_{i=2}^n \int_{t_{i-1}}^{t_i} \Big({\frac{d}{dt}\Big[\frac{\Phi^{(J)}(t)}{v(t)}\Big]}  \Big(\frac{d}{dt}{q(t) }\Big)^{-1}  - \mu_i\Big)\Big(\frac{d}{dt}\Big[\frac{{f(t)}}{v(t)}\Big]\Big)dt\Big)^T  \Big
 \} \ , \nonumber
\end{eqnarray}
and  the constant $k^{(J)}$  is given by
\begin{equation}\label{kJ} \nonumber
k^{(J)} = \begin{cases} \frac{\|\theta^{(J)}\|^4}{(1+c^{(J)})^2}, & \text{ in case (A) }\\
 \frac{\|\theta^{(J)}\|^4}{(1+m^{(J)})^2}, & \text{ in case (B) } . \end{cases}
\end{equation}
 \end{proposition}
 Note that the expected $L^2$-distance in \eqref{reformulated_mse} only differs in the multiplicative factor $k^{(J)}$ for the different cases
 (A) and (B)   and this factor does not depend on the vector-weights $\mu_2, \ldots,\mu_n$.
Therefore optimal weights  minimizing the expected $L^2$-distance can be determined  without distinguishing between the two cases (A) and (B).
 \\
The function $B$ in the criterion \eqref{reformulated_mse} still depends on the unknown regression function $f$ which we replace again by its truncation $f^{(J)}$ defined in \eqref{f_truncated}.  The resulting criterion
is given by
\begin{eqnarray} \label{hol8}
\Phi (\mu_2, \ldots , \mu_n) =
 V (\mu_2, \ldots ,\mu_n)  + B^{(J)}  (\mu_2, \ldots ,\mu_n) ,
\end{eqnarray}
 where
\begin{eqnarray} \label{hol9}
  B^{(J)}  (\mu_2, \ldots ,\mu_n)  &=&
   \mbox{tr}\Big \{
  \sum_{i=2}^n \int_{t_{i-1}}^{t_i} \Big({\frac{d}{dt}\Big[\frac{\Phi^{(J)}(t)}{v(t)}\Big]}  \left(\frac{d}{dt}{q(t) }\right)^{-1} - \mu_i\Big)\Big(\frac{d}{dt}\Big[\frac{{f^{(J)}(t)}}{v(t)}\big]\Big)dt\,   \\
&& \qquad \times \Big(\sum_{i=2}^n \int_{t_{i-1}}^{t_i} \Big({\frac{d}{dt}\Big[\frac{\Phi^{(J)}(t)}{v(t)}\Big]}  \Big(\frac{d}{dt}{q(t) }\Big)^{-1}  - \mu_i\Big)\Big(\frac{d}{dt}\Big[\frac{{f^{(J)}(t)}}{v(t)}\Big]\Big)dt \,  \Big)^T  \Big \}. \nonumber
\end{eqnarray}
We now determine the optimal weights such that the term $ B^{(J)}  (\mu_2, \ldots ,\mu_n) $ in \eqref{hol8} vanishes for all potential Fourier coefficients  $\theta_1, \ldots, \theta_J$ in the function $f^{(J)}$. Therefore, the optimal weights are obtained by minimizing
$\Phi$ in \eqref{hol8}  under the constraint
 \begin{equation}\label{same_expect_trunc}
\int_{0}^{1}\Big[ {\Big(\frac{d}{dt}\Big[\frac{\Phi^{(J)}(t)}{v(t)}\Big]\Big)}  \left(\frac{d}{dt}{q(t) }\right)^{-1}  \Big(\frac{d}{dt}\Big[\frac{{\Phi^{(J)}(t)}}{v(t)}\Big]\Big)^T dt =  \sum_{i=2}^n \mu_i  \int_{t_{i-1}}^{t_i} \frac{d}{dt}\Big[\frac{{\Phi^{(J)}(t)}}{v(t)}\Big]dt.
 \end{equation}
In this situation the criterion \eqref{hol8} reduces  to the minimization of
 \begin{equation}\label{reformulated_mse_trunc_unbiased}
  \mbox{tr}\Big \{\sum_{i=2}^{n}\int_{t_{i-1}}^{t_i} \Big({\frac{d}{dt}\Big[\frac{\Phi^{(J)}(t)}{v(t)}\Big]}   \Big(\frac{d}{dt}{q(t) }\Big)^{-1} - \mu_i\Big) \Big({\frac{d}{dt}\Big[\frac{\Phi^{(J)}(t)}{v(t)}\Big]}   \Big(\frac{d}{dt}{q(t) }\Big)^{-1} - \mu_i\Big)^T \Big(\frac{d}{dt}{q(t) }\Big)dt \Big \},
\end{equation}
with respect to the weights $\mu_2, \ldots, \mu_n$ (depending on the time points $0=t_1<t_2, \ldots, t_{n-1}<t_n=1$). In order to simplify this optimization  we introduce the following notation
 \begin{equation}\label{newnot}
 \beta_i= \frac{\tfrac{\Phi^{(J)}(t_i)}{v(t_i)} - \tfrac{\Phi^{(J)}(t_{i-1})}{v(t_{i-1})}}{\sqrt{q(t_i) - q(t_{i-1})}} \, \quad ,  \quad \gamma_i= \mu_i \sqrt{q(t_i) - q(t_{i-1})} \, 
 \end{equation}
 which however does not reflect the dependence on the time points.
 Using the notation in \eqref{newnot}, the approximation of the expected $L^2$-distance in \eqref{reformulated_mse_trunc_unbiased} can be rewritten in terms of
 the quantities $\gamma_2, \ldots, \gamma_n$ as
 \begin{equation}\label{tr_reformulated_mse_trunc}
\Psi(\gamma_2, \ldots, \gamma_n)  =    -\mbox{tr}(M^{(J)}) +\sum_{i=2}^n { \gamma_i}^{T}{ \gamma_i} ,
 \end{equation}
and the constraint \eqref{same_expect_trunc} is given by
 \begin{equation}\label{extra_cond}
 M^{(J)}= \sum_{i=2}^{n} \gamma_i{\beta_i}^T,
 \end{equation}
where $M^{(J)}$ is  the matrix  defined in \eqref{MJmat}  (for both cases (A) and (B)).
 Note that both the function $\Psi$ and the constraint in \eqref{extra_cond}  do not involve  the function $f$ and only include assumptions concerning the first $J$ basis functions $\varphi_1, \ldots, \varphi_J$ used in the approximation $f^{(J)}$. \\
The resulting optimization problem  \eqref{tr_reformulated_mse_trunc} with constraint \eqref{extra_cond}
 has the same structure as an  optimization problem considered in \cite{dette_new_2017} and from the results in this paper we obtain the solution
\begin{equation}  \label{gammaistar}
\gamma^*_i= M^{(J)} B^{-1} \beta_i, \quad i=2, \ldots, n,
\end{equation}
where the matrix $B$ is given by
\begin{equation}\label{Bmatrix}
B= \sum_{i=2}^n  \beta_i \beta_i^T,
\end{equation}
and $M^{(J)}$ and $\beta_i$ are defined in  \eqref{MJmat} and in \eqref{newnot}, respectively.  If the matrix $B$ is singular, we replace the inverse $B^{-1}$ in \eqref{gammaistar} by a generalized inverse $B^{-}$.
Using the relation between $\gamma_i$ and $\mu_i$ in \eqref{newnot}, we obtain the optimal weights
$$
\mu^*_i = \frac{1}{\sqrt{q(t_i) - q(t_{i-1})}}  M^{(J)} B^{-1} \beta_i, \quad i=2, \ldots n.
$$
Note that these weighs still depend on the design points $t_{2} , \ldots , t_{n-1} $ which will be determined next.

\subsection{Optimal designs for series estimation}  \label{sec43}

Using  the optimal $\gamma^*_2, \ldots, \gamma_n^*$ given in \eqref{gammaistar} in the expression for the function $\Psi$ defined in  \eqref{tr_reformulated_mse_trunc}, we obtain an appropriate
optimal design criterion  for the choice of  the time points $0=t_1 < t_2 < \ldots t_{n-1} < t_n = 1$. More precisely,  for the optimal weights, the function $\Psi$ depends only on the design points and can be represented as the
function
\begin{equation}\label{tilde_psi}
\tilde{\Psi} (t_2, \ldots, t_{n-1})= \mbox{tr}\{M^{(J)}B^{-1}M^{(J)} \},
\end{equation}
 where the matrices $B$  and $M^{(J)}$  are defined in  \eqref{Bmatrix} and \eqref{MJmat}, respectively and
depend on $0=t_1< t_2, \ldots, t_{n-1}< t_n=1$. The optimal design is now determined by minimizing the function $\tilde \Psi$, which
 is different from the criterion considered  in \cite{dette_new_2017} for unbiased linear estimation in the linear
regression model
\begin{equation}\label{auxil_mod}
Y_{t_i} = (\Phi^{(J)}(t_i))^T \theta + \varepsilon_{t_i}, \quad i=1, \ldots, n .
\end{equation}
The optimal time points only depend on the first $J$ basis functions which are used for the estimator of the  regression function $f$ and have to be determined numerically in all cases of practical interest. We will present some examples in Section
\ref{sec5}.

\subsection{The final estimate}  \label{sec44}

With the optimal weights $\mu_2^*, \ldots , \mu_n^*$  determined in Section \ref{sec42} and the optimal time points $t_2^*, \ldots t_{n-1}^*$ determined in Section \ref{sec43}, the estimators  in \eqref{hol10} and \eqref{hol11}
corresponding to the cases (A) and (B) are given by
\begin{equation}\label{fin_lin_estA}
 \hat\theta^{(J), n}=
\frac{1}{1+c^{(J)}}\theta^{(J)}(\theta^{(J)})^{T} \big\{M^{(J)} B^{-1} \sum_{i=2}^{n}  \beta_i^{*}
(\frac{Y_{t^*_i}}{v(t^*_i)} - \frac{Y_{t^*_{i-1}}}{v(t^*_{i-1})} ) + \frac{\Phi^{(J)}(0)}{u(0)} \frac{Y_0}{v(0)} \big\},    \end{equation}
 and
 \begin{equation}\label{fin_lin_estB}
 \hat\theta^{(J), n}=
\frac{1}{1+m^{(J)}}\theta^{(J)}(\theta^{(J)})^{T} M^{(J)} B^{-1}\sum_{i=2}^{n} \beta_i^{*}
 (\frac{Y_{t^*_i}}{v(t^*_i)} - \frac{Y_{t^*_{i-1}}}{v(t^*_{i-1})} ) ~,
 \end{equation}
 respectively,
 where
$$
   \beta^*_i= \frac{\tfrac{\Phi^{(J)}(t^*_i)}{v(t^*_i)} - \tfrac{\Phi^{(J)}(t^*_{i-1})}{v(t^*_{i-1})}}{\sqrt{q(t^*_i) - q(t^*_{i-1})}}  \quad i=2, \ldots n.
 $$

 For their application we still require knowledge of the vector of Fourier coefficients $\theta^{(J)}$ and the constants $c^{(J)}$ and $m^{(J)}$ defined in \eqref{cJ2} and \eqref{cJ2u0} (note that these quantities
 also depend on $\theta^{(J)}$).
For this purpose we  propose to use the linear unbiased estimate derived by \cite{dette_new_2017} for the linear model \eqref{auxil_mod}.   This estimate is defined as
  \begin{equation}\label{blue_approx}
 \check\theta^{(J),n} =
 (C^{(J)})^{-1}\Big\{ M^{(J)} B^{-1} \sum_{i=2}^{n} \beta_i^{*}  (\frac{Y_{t^*_i}}{v(t^*_i)} - \frac{Y_{t^*_{i-1}}}{v(t^*_{i-1})} ) + \frac{\Phi^{(J)}(0)}{u(0)} \frac{Y_0}{v(0)} \Big\},
  \end{equation}
  and the quantity $c^{(J)}$ in \eqref{cJ2}  is estimated by  $\check c^{(J), n}  = (\check\theta^{(J),n})^{T}C^{(J)} \check\theta^{(J),n}   $. A straightforward calculation shows  that the resulting estimator for the case (A)
   is given by
\begin{eqnarray}  \nonumber
 \hat\theta^{(J), n} &=& \frac{1}{1+\check{c}^{(J),n}}\check\theta^{(J),n}(\check\theta^{(J),n})^{T}\big\{ M^{(J)} B^{-1} \sum_{i=2}^{n} \beta_i^{*}  (\frac{Y_{t^*_i}}{v(t^*_i)} - \frac{Y_{t^*_{i-1}}}{v(t^*_{i-1})} ) + \frac{\Phi^{(J)}(0)}{u(0)} \frac{Y_0}{v(0)} \big\}\\
 &= &   \frac{1}{1+\check{c}^{(J),n}}\check\theta^{(J),n}(\check\theta^{(J),n})^{T}C^{(J)} \check\theta^{(J),n} =   \frac{\check{c}^{(J), n}}{1+\check{c}^{(J), n}}\check\theta^{(J),n},\label{final_lin_est3}
 \end{eqnarray}
which is a shrinkage version of the estimator $\check{\theta}^{(J), n}$ in \eqref{blue_approx}.\\
For the case (B) similar arguments show that the estimator in \eqref{fin_lin_estB} can also be rewritten in terms of the linear unbiased estimate $\check\theta^{(J),n}$, that is,
\begin{eqnarray*}
\hat\theta^{(J), n}&=&
\frac{1}{1+\check{m}^{(J)}}\check\theta^{(J)}(\check\theta^{(J)})^{T} M^{(J)} B^{-1}\sum_{i=2}^{n} \beta_i^{*}  (\frac{Y_{t^*_i}}{v(t^*_i)} - \frac{Y_{t^*_{i-1}}}{v(t^*_{i-1})} )
=  \frac{\check{m}^{(J), n}}{1+\check{m}^{(J), n}}\check\theta^{(J),n},
\end{eqnarray*}
where $\check m^{(J), n}  = (\check\theta^{(J),n})^{T}M^{(J)} \check\theta^{(J),n}$. Here the structure of the estimator $\check{\theta}^{(J), n}$ depends on the structure of the basis functions contained in the vector $\Phi^{(J)}$ [see Section 5 in \cite{dette_new_2017} for more details].

\section{Numerical results}\label{sec5}
\label{sec5}
\def\theequation{5.\arabic{equation}}
\setcounter{equation}{0}

In this  section, we illustrate the properties of the estimator and the corresponding optimal design derived in Section \ref{sec4} by means of a small simulation study. We consider a Gaussian process assuming both an exponential kernel and a Brownian motion as the error process in model \eqref{disc_mod}.
In both cases, we present the numerically calculated optimal time points with respect to the criterion defined in \eqref{tilde_psi} and the corresponding simulated integrated mean squared errors for the estimator
\begin{equation}\label{es1}
\hat f^{(J), n} (t)= \sum_{j=1}^{J}\hat \theta^{(J),n} \varphi_j(t),  \\
\end{equation}
proposed in this paper and the estimator
\begin{equation}\label{es2}
\check f^{(J), n}(t)  = \sum_{j=1}^{J} \check\theta^{(J),n} \varphi_j(t) ,
\end{equation}
which is based on the best linear unbiased estimates in the tuncated Fourier expansion.  \\
Throughout this section, we will use the trigonometric series defined in \eqref{cosine_basis} as orthonormal basis of $L^2( [0, 1]) $. We further assume that the unknown function $f$ is symmetric on the interval $[0, 1]$ such that it is sufficient to use only  the cosine functions  in the series expansions  of $f$.  Consequently, the orthonormal system is given by
$$
\varphi_1(t) = 1 \ , \quad \varphi_{j}(t) =  \sqrt{2}\cos(2\pi (j-1) t) \ ,  \quad j=2, 3, \ldots \
$$
In Section \ref{num_exp} we consider the exponential kernel, whereas in Section \ref{num_brown} we concentrate on the Brownian motion.

\subsection{The Exponential kernel}\label{num_exp}
We assume that the error process $\{\varepsilon_t: t \in[0, 1]\}$ is a centered Gaussian process with an exponential kernel of the form
$K(s,t) = \exp(-L|s-t|)$, where $L\in \mathbb{R}^+$ is a given constant. This  can be represented in the triangular form \eqref{triangular-kernel}
with $u(t) =\exp(Lt)$ and $v(t) = \exp(-Lt)$    and
the function $q$ is obtained as $q(t) = u(t)/v(t) = \exp(2Lt)$.
Therefore, we have  $u(0)\neq0$ (which  corresponds to case (A)) and the preliminary estimator $\check\theta^{(J), n}$ in \eqref{blue_approx}   is given by
 \begin{equation*}
 \check\theta^{(J),n} = (C^{(J)})^{-1} \big\{M^{(J)} B^{-1}\sum_{i=2}^n \frac{e^{Lt_i} \Phi^{(J)}(t_i^*) - e^{Lt_{i-1}} \Phi^{(J)}(t^*_{i-1})}{\sqrt{e^{2Lt_i} - e^{2Lt_{i-1}} }}\Big(e^{Lt_i} Y_{t^*_i} - e^{Lt_{i-1}} Y_{t^*_i-1} \Big)+ {Y_0}  \Phi^{(J)}(0) \big\},
 \end{equation*}
where the matrices $M^{(J)}$, $C^{(J)}$ and $B$ become
\begin{eqnarray}
M^{(J)}& =& \int_0^1 \frac{ (\dot{\Phi}^{(J)}(t) +L{\Phi}^{(J)}(t))(\dot{\Phi}^{(J)}(t) + L{\Phi}^{(J)}(t))^T} {2L} dt, \nonumber \\
C^{(J)} &= &M^{(J)}+ \Phi^{(J)}(0)(\Phi^{(J)}(0))^T,\nonumber \\
B  &= & \sum_{i=2}^{n} \frac{(e^{Lt_i} \Phi^{(J)}(t_i) - e^{Lt_{i-1}} \Phi^{(J)}(t_{i-1})) (e^{Lt_i} \Phi^{(J)}(t_i) - e^{Lt_{i-1}} \Phi^{(J)}(t_{i-1}))^T}{e^{2Lt_i} - e^{2Lt_{i-1}}}.\nonumber \label{BLmatrix}
\end{eqnarray}
The  estimator $\hat\theta^{(J), n}$ proposed in this paper is given in equation \eqref{final_lin_est3} and the corresponding estimators of the function $f$ are defined in \eqref{es1} and \eqref{es2}.

We first consider the  exponential covariance kernel  with $L=1$
and assume that three basis functions  $\varphi_1(t) = 1,\ \varphi_2(t)=\cos(2\pi t),\ \varphi_3(t) = \cos(4�\pi t)$
are used in the series estimator, where   $n=4$ and $n=7$ observations at different
 time points $0=t_1 < t_2 <  \ldots <  t_{n-1} <t_n=1$ can be taken. Note that one needs at least $n=4$ observations
 at different time points to guarantee that the matrix  $B$ in the preliminary estimator $ \check\theta^{(J),n}$   is non-singular. The optimal points are determined minimizing the  criterion  \eqref{tilde_psi}
 by particle swarm optimization [see \cite{clerc_particle_2006} for details] and the results
are presented in the first row of Table \ref{tab_tp}. The second row shows the results for $L=5$ and interestingly the optimal points do not change substantially for different values of the constant $L$. Also note that all designs are nearly equidistant.

\vspace{0.5cm}
\begin{table}[h!]
\centering
 \renewcommand{\arraystretch}{1.25}
\begin{tabular}{| c | c | c|}
  \hline
$L$  & $n=4$ & $n=7$ \\
  \hline
$1$ &{$0.00, \, 0.25, \, 0.52, \, 1.00$ }	&{$0.00  , \, 0.12  , \, 0.27  , \, 0.45  , \, 0.57 , \,   0.77  , \,  1.00$}\\
 \hline
 $5$ &{$0.00, \, 0.25, \, 0.51, \, 1.00$} &{$0.00 , \, 0.12 , \, 0.27, \,  0.45, \,  0.57, \,  0.76, \,  1.00 $}\\
\hline
\end{tabular}
\caption{\it Optimal time points for series estimation  minimizing the criterion \eqref{tilde_psi}. The covariance kernel is given by   $\exp(-L|s-t|)$ with $L=1$
(first row) and $L=5$ (second row).
 \label{tab_tp}}
\end{table}
%
%

We now evaluate the performance of the different estimators and the optimal time points
by means of a simulation study. For the sake of comparison we also consider  non-optimized
  time points for the simulation, which are given by
  \begin{eqnarray} \label{comp4}
  &&0.00 , \, 0.45 , \,  0.90 , \,  1.00  \\
\label{comp7}
&& 0.00 , \,0.18, \, 0.36 , \,0.54, \, 0.72 , \,0.90 , \,1.00
  \end{eqnarray}
  for the case $n=4$ and $n=7$, respectively.

In the simulation study  we generate data  according to model \eqref{disc_mod}  with two regression functions
  \begin{eqnarray}
   \label{simmod1} f(t) &=& 4t(t-1),  \\
 \label{simmod2}  f(t) &=& \sqrt{t(t-1)},
  \end{eqnarray}
(note that both proposed functions are symmetric with $f(0)=f(1)=0$).
For each model
the mean integrated squared error of  the estimators $\hat f^{(J), n}$ and $\check f^{(J), n}$
defined in \eqref{es1} and \eqref{es2} respectively is  determined.
More precisely, if $S$ denotes the number of simulation runs and
$\bar f_\ell$ is the estimator based on the $\ell$-th run (either  $\hat f^{(J), n}$ and $\check f^{(J), n}$), the simulated mean integrated squared
error, $\text{MISE}_n$, is given by
$$\text{MISE}_n= \frac{1}{S}\sum_{\ell=1}^{S} \int_0^1 \big(\bar f_{\ell}(t) - f(t) \big)^2 dt, $$
where $f$, the ``true''  regression  function  under consideration, is either given by \eqref{simmod1}  or by  \eqref{simmod2}.
All results are based on $S=1000$ simulation runs.


\vspace{0.5cm}
\begin{table}[h!]
\centering
 \renewcommand{\arraystretch}{1.25}
 \begin{tabular}{| c | c | c| c ||c |c | }
 \hline
 \multicolumn{2}{|c|}{ }  &\multicolumn{2}{c||}{$n=4$} &\multicolumn{2}{c|}{$n=7$}  \\
  \hline
\multicolumn{2}{|c|}{ }  &\multicolumn{2}{c||}{design} &\multicolumn{2}{c|}{design}  \\
  \hline
 \hline
$f$ & estimator  & optimal & \eqref{comp4} & optimal & \eqref{comp7} \\
  \hline
\multirow{2}{*}{\eqref{simmod1}  } & $\hat{f}^{(J), n}$ & 1.72 & 2.06 & 1.58 & 1.59 \\
\cline{2-6}
&  $\check{f}^{(J), n}$  & 1.89 & 2.22 & 1.76 & 1.77\\
     \hline
     \hline
\multirow{2}{*}{\eqref{simmod2} } & $\hat{f}^{(J), n}$ & 1.67 & 2.04 & 1.54 & 1.56 \\
\cline{2-6}
& $\check{f}^{(J), n}$  & 1.89 & 2.21 & 1.76 & 1.79 \\
\hline
\end{tabular}
\caption{\it Simulated mean integrated squared error  of the estimators $\hat{f}^{(J), n}$   and $\check{f}^{(J), n}$
defined in \eqref{es1}  and \eqref{es2} for different regression functions. The covariance kernel is given by $\exp(-|s-t|)$. Third column:  optimal design; Fourth column:  comparative design in \eqref{comp4}.  Left part: $n=4$ observations; right part:  $n=7$ observations. 
 \label{tab_sim_exp1-n4}}
\end{table}

For the case of the sample size $n=4$, the resulting mean integrated squared error  of the different estimators (and corresponding optimal time points)
is shown in the left part of Table \ref{tab_sim_exp1-n4}. For instance, the mean integrated squared error   of the estimator  $\hat f^{(J), n}$ (based on the
on the optimal design)  is   $1.72$, if the  true function is given by  \eqref{simmod1}, whereas it is
 $2.06$ if the observations are taken according to the non-optimized design \eqref{comp4}.
 Thus, the optimal design yields a reduction by
 $17\%$ in the mean integrated squared error.
 The optimal design also yields a reduction  of $15\%$ of  the mean squared error of the preliminary estimator $\check f^{(J), n}$ (although it is not constructed for this purpose).
 We also observe that the new  estimator $\hat f^{(J), n}$ clearly outperforms the estimator $\check f^{(J), n}$ in
 all cases   under consideration (reduction of the mean squared error between $9\%$ and $12\%$).

For the case of the sample size $n=7$, the corresponding results  are presented in the right part of Table \ref{tab_sim_exp1-n4} and  we observe a similar behavior.
The new  estimator  $\hat f^{(J), n}$  clearly outperforms  $\check f^{(J), n}$ regardless of the design and model  under consideration.
On the other hand the improvement by the choice of the design is less visible compared to the
case where the sample size is $n=4$. This means that the influence of design on the performance of the estimators decreases with increasing sample size.
The reason for this observation lies in the fact that in the models under consideration   the discrete model  \eqref{disc_mod} already
provides a good approximation of the continuous model \eqref{cont_mod} for the sample size $n=7$. As in this   model the full trajectory is available
  the impact of the design is   negligible for sample sizes larger than $10$. As a consequence, a larger sample size would not decrease the integrated mean squared error
  substantially either.
A similar effect was also observed by \cite{dette_new_2017} in the linear regression model with correlated observations.

\vspace{0.5cm}
\begin{table}[h!]
\centering
 \renewcommand{\arraystretch}{1.25}
 \begin{tabular}{| c | c | c| c ||c |c | }
 \hline
 \multicolumn{2}{|c|}{ }  &\multicolumn{2}{c||}{$n=4$} &\multicolumn{2}{c|}{$n=7$}  \\
  \hline
\multicolumn{2}{|c|}{ }  &\multicolumn{2}{c||}{design} &\multicolumn{2}{c|}{design}  \\
  \hline
   \hline
 $f$  & estimator  & optimal & \eqref{comp4} & optimal & \eqref{comp7}\\
  \hline
\multirow{2}{*}{\eqref{simmod1}  }  & $\hat{f}^{(J),n}$  & 0.65 & 2.13 & 0.47 & 0.51\\
\cline{2-6}
&  $\check{f}^{(J),n}$  &0.77 & 2.30 & 0.58 & 0.62\\
     \hline
     \hline
\multirow{2}{*}{\eqref{simmod2} }  & $\hat{f}^{(J),n}$   &0.64 & 2.09 & 0.43 & 0.43 \\
\cline{2-6}
& $\check{f}^{(J),n}$  & 0.81 & 2.30 & 0.59 & 0.59 \\
\hline
\end{tabular}
\caption{\it Simulated mean integrated squared error  of the estimators $\hat{f}^{(J), n}$  and $\check{f}^{(J), n}$
defined in \eqref{es1}  and \eqref{es2} for different regression functions. The covariance kernel is given by $\exp(-5|s-t|)$. Third column:  optimal design; Fourth column:  comparative  design in \eqref{comp4}. Left part: $n=4$ observations; right part: $n=7$ observations. \label{tab_sim_exp5-n4}}
\end{table}

Next we consider a situation where  the correlation between the different observations is smaller and so we use  the constant $L=5$ for the exponential kernel.
The time points minimizing the criterion $\tilde\Psi$ in \eqref{tilde_psi} are depicted in the second row of Table \ref{tab_tp} for $n=4$ and $n=7$.
The simulated mean integrated squared error   of the estimators $\hat f^{(J), n}$ and $\check f^{(J), n}$
defined in \eqref{es1} and \eqref{es2}  are displayed in Table \ref{tab_sim_exp5-n4}  for the cases of sample size $n=4$ and $n=7$.
When the sample size is $n=4$, we observe that the optimal design yields a substantial reduction  in the mean squared errors of both estimators (between $65\%$ and $70\%$).
Compared to the case $L=1$ (see Table \ref{tab_sim_exp1-n4}) the  reduction is  larger.
When the   sample size is $n=7$ the mean integrated squared error of the estimators based on the optimal time points are slightly smaller compared to the  non-optimized time points. We observe again that the influence of the position of the time points, and thus of the design, decreases if the sample size $n$ increases (see Table \ref{tab_sim_exp5-n4}).
A comparison of  the two  estimators  \eqref{es1}  and \eqref{es2}  shows again that the new estimator $\hat{f}^{(J), n}$ outperforms the estimator $\check f^{(J), n}$ in all cases under consideration (reduction of the mean squared error between $16\%$ and $27\%$).

%

\subsection{Brownian motion}\label{num_brown}
We now consider the case where the error process in \eqref{disc_mod} is given by a Brownian motion, that is
$K(s,t) = s \wedge t \,$, which can be represented by
$
K(s,t)   = s, \quad s \leq t .
$
Therefore,  the functions $u$ and $v$  in \eqref{triangular-kernel} are given by $u(t) =t$ and $v(t) = 1$, respectively, and the function $q$ is obtained as $q(t)=u(t)/v(t)= t$.
This situation corresponds to case (B), where $u(0) = 0 $ and $f(0) = 0$.
The estimator $\hat{\theta}^{(J), n} $ is given by
\begin{equation}\label{linest_brown}
\hat{\theta}^{(J), n} = \frac{1}{1+\check m^{(J)}}\check\theta^{(J)}(\check\theta^{(J)})^{T} M^{(J)} B^{-}\sum_{i=2}^{n}\frac{(\Phi^{(J)}(t_i) - \Phi^{(J)}(t_{i-1}))^T}{\sqrt{{t_i} - t_{i-1}}} ({Y_{t_i}}- Y_{t_{i-1} } ),
\end{equation}
where the matrices $M^{(J)}$, $B$  and the constant $m^{(J)}$ are of the form
\begin{eqnarray*}
M^{(J)}& =& \int_0^1 {\dot\Phi^{(J)}}(t)({\dot\Phi^{(J)}}(t))^T dt, \\
B  &= & \sum_{i=2}^{n} \frac{(\Phi^{(J)}(t_i) -\Phi^{(J)}(t_{i-1})) (\Phi^{(J)}(t_i) - \Phi^{(J)}(t_{i-1}))^T}{{t_i} - t_{i-1}} ,\\
\check{m}^{(J)} &=& (\check\theta^{(J)})^T M^{(J)} \,  \check\theta^{(J)}.
\end{eqnarray*}
Note that both the first row and the first column of the matrices $M^{(J)}$ and $B$ are zero (since $\varphi_1(t) =1$), such that both matrices are singular. Consequently,  as proposed in Section \ref{sec4},
we use
the generalized inverse
$$B^- = \begin{pmatrix} 0 & 0 \\ 0 & \tilde{B}^{-1} \end{pmatrix}, $$
 of $B$, where the matrix $\tilde{B}$  is given by
$$\tilde{B} = \begin{pmatrix}\bf{0}_{(J-1)} & \mathbf{I}_{(J-1)\times (J-1)} \end{pmatrix} B \begin{pmatrix}\bf{0}^T_{(J-1)} \\ \mathbf{I}_{(J-1)\times (J-1)} \end{pmatrix} .
$$ 
Here the vector $\bf{0}_{(J-1)}$ is of dimension  $(J-1)$  with zero entries.
where the matrix $\mathbf{I}_{(J-1)\times (J-1)}$ is the $(J-1)$ dimensional identity matrix. 
The estimator $\check\theta^{(J),n}$ is obtained from Section 5.2 in \cite{dette_new_2017}
\begin{eqnarray*}
\check\theta^{(J),n} = \underline{C}^{(J)} \sum_{i=2}^{n}\frac{(\Phi^{(J)}(t_i) - \Phi^{(J)}(t_{i-1}))^T}{\sqrt{{t_i} - t_{i-1}}} ({Y_{t_i}}- Y_{t_{i-1} } ),\end{eqnarray*}
where the matrix $ \underline{C}^{(J)} $ is of the form
$$
 \underline{C}^{(J)} =
 \begin{pmatrix}
 0 & -(\Phi^{(J)}(0))^T \begin{pmatrix}0 \\�\tilde{B}^{-1} \end{pmatrix} \\
\bf{0}_{(J-1)}  & \tilde{B}^{-1}
 \end{pmatrix}.
$$  

We now analyze the behavior of the resulting estimators of the function $f$ if the first three basis functions are used for the series estimator and $n=4$ or $n=7$ observations at different time points $0= t_1<t_2 < \ldots < t_{n-1} < t_n = 1$ are available. The optimal time points minimizing the criterion \eqref{tilde_psi} derived in Section \ref{sec4} are given by
\begin{eqnarray}\label{optpoints_brownian_n4}
&& 0.00,\, 0.25,\, 0.47, \, 1.00 \\
&&
\label{optpoints_brownian_n7}
0.00 , \, 0.22 , \,  0.28  , \, 0.50 , \,  0.72 , \,  0.78 , \, 1.00
\end{eqnarray}
for sample sizes $n=4$ and $n=7$ respectively.
Note that the optimal time points  \eqref{optpoints_brownian_n4} and  \eqref{optpoints_brownian_n7} differ from the optimal time points for the case of the exponential kernel displayed in Table \ref{tab_tp}.  This indicates that  the position of the optimal time points depends on the structure of the covariance kernel.  \\

\begin{table}[h!]
\centering
 \renewcommand{\arraystretch}{1.25}
 \begin{tabular}{| c | c | c| c ||c |c | }
 \hline
 \multicolumn{2}{|c|}{ }  &\multicolumn{2}{c||}{$n=4$} &\multicolumn{2}{c|}{$n=7$}  \\
  \hline
\multicolumn{2}{|c|}{ }  &\multicolumn{2}{c||}{design} &\multicolumn{2}{c|}{design}  \\
  \hline
   \hline
 $f$  & estimator  & optimal & \eqref{comp4} & optimal & \eqref{comp7}\\
  \hline
\multirow{2}{*}{\eqref{simmod1}  }  & $\hat{f}^{(J),n}$  & 0.16 & 0.41 & 0.13 & 0.14 \\
\cline{2-6}
&  $\check{f}^{(J),n}$  & 0.15 & 0.43 & 0.12 & 0.12 \\
     \hline
     \hline
\multirow{2}{*}{\eqref{simmod2}  }  & $\hat{f}^{(J),n}$  & 0.13 & 0.45 & 0.11 & 0.11 \\
\cline{2-6}
& $\check{f}^{(J),n}$  & 0.15 & 0.48  & 0.12 & 0.13 \\
\hline
\end{tabular}
\caption{\it Simulated mean integrated squared error  of the estimators $\hat{f}^{(J), n}$  and $\check{f}^{(J), n}$
defined in \eqref{es1}  and \eqref{es2} for different regression functions. The error process is a Brownian motion. Third column:  optimal design; Fourth column:  comparative  design in \eqref{comp4}. Left part: $n=4$ observations; right part:  $n=7$ observations.   }  \label{tab_sim_brown-n4}
\end{table}


 The resulting mean integrated squared errors  of  the estimators $\hat f^{(J), n}$ and $\check f^{(J), n}$ are displayed in Table \ref{tab_sim_brown-n4}, where we again consider the comparative set of time points depicted in \eqref{comp4} and \eqref{comp7}. We obtain similar results as in Section \ref{num_exp}. More specifically, for the case of sample size $n=4$, we observe that the optimal design yields a substantial reduction  in the mean squared errors of both estimators (see the left part of Table \ref{tab_sim_brown-n4}).
When the sample size is $n=7$, the difference between the optimal time points and the design \eqref{comp7} is less visible.  \\
A comparison of the two estimators  shows  a different behavior as in Section \ref{num_exp}, that is, unlike the case of an exponential Kernel, when the error process is a Brownian motion, both estimators perform well and they have similar (small) mean integrated squared errors (see Table \ref{tab_sim_brown-n4}). \\

\bigskip

{\bf Acknowledgements} This work has been supported in part by the
Collaborative Research Center ``Statistical modeling of nonlinear
dynamic processes'' (SFB 823, Teilprojekt C2) of the German Research Foundation
(DFG) and  by a grant from the National Institute of General Medical Sciences of the National
Institutes of Health under Award Number R01GM107639. The content is solely the responsibility of the authors and does not necessarily
 represent the official views of the National
Institutes of Health.

\begin{small}
 \bibliography{nonparambib}
\end{small}

\section{Technical details} \label{sec6}
{\bf Proof of Theorem \ref{theo_cont_opt_est}}
We restrict ourselves to the proof of the result in case (A), the other cases can be proved in a similar way.
Note that the function $\Psi_j$ is convex on the space of all signed measures and therefore, a signed measure  $\xi_j^*$ minimizes $\Psi_j$ if and only if
the directional derivative from $\xi^*_j$ in any direction is nonnegative, that is
$$\frac{\partial}{\partial \alpha} \Psi_j((1-\alpha)\xi^*_j + \alpha \eta ) \Big |_{\alpha=0} \geq   0,$$
for all signed measures $\eta$ on the interval $[0,1]$.
A straightforward calculation gives
\begin{equation}\label{deriv}
\begin{split}
\frac{\partial}{\partial \alpha} \Psi_j((1-\alpha)\xi^*_j + \alpha \eta ) \Big |_{\alpha=0} =& \int_{0}^{1} \int_{0}^1 f(s) f(t) + K(s,t) \big(\xi^*_j(ds) \xi^*_j(dt) - \xi^*_j(dt) \eta(dt)\big) \\
& + \theta_j \int_0^1 f(t) \big(\eta(dt) - \xi^*_j(dt)\big).
\end{split}
\end{equation}
Consequently, the signed measure $\xi^{*}_j$ minimizes $\Psi_j$ if and only the inequality
\begin{equation}\label{equiv_ineq}
\int_{0}^{1} \int_{0}^1 f(s) f(t) + K(s,t) \big( \xi^*_j(ds) \xi_j(dt) -  \xi^*_j(dt) \eta(dt)\big) + \theta_j \int_0^1 f(t) \big(\eta(dt) -  \xi^*_j(dt)\big) \geq  0,
\end{equation}
is satisfied for all signed measures $\eta$ on the interval $[0,1]$. \\
In order to check \eqref{equiv_ineq} for the signed measure $\xi^*_j$ we calculate each term in \eqref{equiv_ineq} separately, where we use the following representation of the quantities in \eqref{c}-\eqref{pt}
\begin{eqnarray*}
c&=& \int_{0}^{1}  \frac{1}{v^2(t)} \frac{ (\dot{f}(t) v(t)-f(t)\dot{v}(t))^2}{\dot{u}(t)v(t) - u(t)\dot{v}(t)} dt + \frac{f^2(0)}{u(0)v(0)} ,\\
P_0 &=& \frac{1}{u(0)} \frac{f(0)\dot{u}(0) - \dot{f}(0) u(0)}{v(0)\dot{u}(0) - \dot{v}(0) u(0)}, \\
P_1&=& \frac{1}{v(1)}\frac{ \dot{f}(1) v(1)-f(1)\dot{v}(1)}{\dot{u}(1)v(1) - u(1)\dot{v}(1) } ,\\
p(t) &=&- \frac{1}{v(t)}\frac{d}{dt} \left[ \frac{ \dot{f}(t) v(t)-f(t)\dot{v}(t)}{\dot{u}(t)v(t) - u(t)\dot{v}(t) } \right].
\end{eqnarray*}
To simplify \eqref{equiv_ineq} we note that integration by parts yields
\begin{eqnarray*}
\int_{0}^{1}f(s) \xi^*_j(ds) &=&
 \frac{\theta_j}{1+c} \Big(\frac{f(0)(\dot{u}(0) f(0)
  -  u(0)\dot{f}(0))}{u(0)(\dot{u}(0) v(0) - u(0)\dot{v}(0) )} + \frac{f(1)(v(1)\dot{f}(1)-\dot{v}(1) f(1)) }{v(1)(\dot{u}(1) v(1) - u(1)\dot{v}(1) )} \\
&&- \int_0^1 \frac{f(s)}{v(s)} \frac{d}{ds} \left[ \frac{v(s)\dot{f}(s)-\dot{v}(s) f(s)}{\dot{u}(s) v(s) - u(s)\dot{v}(s)} \right]  ds \Big) \\
&=& \frac{\theta_j}{1+c} \Big(\frac{f(0)(\dot{u}(0) f(0)  -  u(0)\dot{f}(0))}{u(0)(\dot{u}(0) v(0) - u(0)\dot{v}(0) )} + \frac{f(1)(v(1)\dot{f}(1)-\dot{v}(1) f(1)) }{v(1)(\dot{u}(1) v(1) - u(1)\dot{v}(1) )}\\
&&- \Big[ \frac{f(s)}{v(s)}\frac{v(s)\dot{f}(s)-\dot{v}(s) f(s)}{\dot{u}(s) v(s) - u(s)\dot{v}(s)}  \Big]^1_0 + \int_0^1 \frac{1}{v^2(s)} \frac{ \big(\dot{v}(s) f(s)  -  v(s)\dot{f}(s)\big)^2 }{\dot{u}(s) v(s) - u(s)\dot{v}(s)}ds\Big) \\
&=& \frac{\theta_j}{1 + c} \Big(\frac{f^2(0)}{u(0) v(0)} + \int_0^1 \frac{1}{v^2(s)} \frac{ \big(\dot{v}(s) f(s)  -  v(s)\dot{f}(s)\big)^2 }{\dot{u}(s) v(s) - u(s)\dot{v}(s)}ds  \Big) \\
&=& \frac{\theta_j}{1+c} c \, .
\end{eqnarray*}
Similarly, we obtain
\begin{eqnarray*}
\int_0^1 K(s, t) \xi^*_j(ds) &=& \int_0^t u(s) v(t) \xi^*_j(ds) +  \int_t^1 u(t) v(s) \xi^*_j(dt)\\
&=& \frac{\theta_j}{1+c} \Big(v(t) \frac{\dot{u}(0) f(0)  -  u(0)\dot{f}(0)}{\dot{u}(0) v(0) - u(0)\dot{v}(0)} - v(t) \int_0^t \frac{u(s)}{v(s) }  \frac{d}{ds} \left[ \frac{v(s)\dot{f}(s)-\dot{v}(s) f(s)}{\dot{u}(s) v(s) - u(s)\dot{v}(s)} \right] ds  \\
&&  + u(t)  \frac{v(1)\dot{f}(1)- \dot{v}(1) f(1)}{\dot{u}(1) v(1) - u(1)\dot{v}(1)}
- u(t)  \int_t^1 \frac{d}{ds} \left[ \frac{v(s)\dot{f}(s)-\dot{v}(s) f(s)}{\dot{u}(s) v(s) - u(s)\dot{v}(s)} \right] ds \Big) \\
&=& \frac{\theta_j}{1+c} \Big(v(t) \frac{f(0)}{v(0)} + f(t) - v(t) \frac{f(0)}{v(0)} \Big) \\
&=& \frac{\theta_j}{1+c} f(t) ,
\end{eqnarray*}
where we have used again integration by parts for the third equality.
Consequently, we get
\begin{eqnarray*}
\int_0^1 \int_0^1 f(s) f(t) \xi^*_j(ds) \xi^*_j(dt) = \Big(\int_{0}^{1}f(s) \xi^*_j(ds) \Big)^2 = \frac{\theta^2_j}{(1+c)^2} c^2,
\end{eqnarray*}
\begin{eqnarray*}
\int_0^1 K(s, t) \xi^*_j(ds) \xi^*_j(dt) = \frac{\theta_j}{1+c}  \int_0^1 f(t) \xi^*_j(dt) =  \frac{\theta^2_j}{(1+c)^2}c ,
\end{eqnarray*}
and thus the left hand side of  \eqref{equiv_ineq} reduces to
$$\frac{\theta^2_jc^2}{(1+c)^2} + \frac{\theta^2_jc}{(1+c)^2} - \frac{\theta^2_j c}{1+c} - \int_{0}^1 f(t) \eta(dt) \big(\frac{\theta_j}{1+c}c + \frac{\theta_j}{1+c} - {\theta_j}\big)= 0,$$
for an arbitrary signed measure $\eta$. This proves that \eqref{equiv_ineq} holds and the signed measure $\xi^*_j$ defined in Theorem \ref{theo_cont_opt_est} minimizes the function $\Psi_j$.

\medskip
{\bf Proof of Proposition \ref{lem_reformulated_mse}}  
For the term on the left hand side of equation \eqref{reformulated_mse} we obtain
\begin{eqnarray}
\mathbb{E}\big[\|\hat\theta^{(J), *} - \hat\theta^{(J), n} \|^2\big]
&=& \mbox{tr}\Big\{  \mathbb{E}\Big[\big(\hat\theta^{(J), *} - \hat\theta^{(J), n} \big) \big(\hat\theta^{(J), *} - \hat\theta^{(J), n} \big)^T\Big] \Big\}  \nonumber \\
&=& \frac{\|\theta^{(J)}\|^4}{(1+c^{(J)})^2}\mbox{tr}\Big\{  \mathbb{E}\Big[\Big(\sum_{i=2}^{n} \int_{t_{i-1}}^{t_i} \Big(\frac{d}{dt} \left[ \frac{\Phi^{(J)}(t)}{v(t)} \right] \big(\frac{d}{dt}q(t) \big)^{-1} - \mu_i \Big)d\Big(\frac{Y_t}{v(t)} \Big) \Big)  \nonumber \\
& & \quad \times
\Big(\sum_{i=2}^{n} \int_{t_{i-1}}^{t_i} \Big(\frac{d}{dt} \left[ \frac{\Phi^{(J)}(t)}{v(t)} \right] \big(\frac{d}{dt}q(t) \big)^{-1} - \mu_i \Big)d\Big(\frac{Y_t}{v(t)} \Big) \Big)^T\Big]\Big\} \nonumber .
\end{eqnarray}
For the determination of the expected value inside the trace
\begin{equation}\label{expec_value}
\mathbb{E}\Big[\Big(\sum_{i=2}^{n} \int_{t_{i-1}}^{t_i} \Big(\frac{d}{dt} \frac{\Phi^{(J)}(t)}{v(t)} \big(\frac{d}{dt}q(t) \big)^{-1} - \mu_i \Big)d\Big(\frac{Y_t}{v(t)} \Big)\Big)
\Big(\sum_{i=2}^{n} \int_{t_{i-1}}^{t_i} \Big(\frac{d}{dt} \frac{\Phi^{(J)}(t)}{v(t)} \big(\frac{d}{dt}q(t) \big)^{-1} - \mu_i \Big)d\Big(\frac{Y_t}{v(t)} \Big) \Big)^T\Big],
\end{equation}
we use a transformation of the Gaussian process $\{Y_t: t\in [0, 1]\}$ to a Brownian motion, as it was introduced by \cite{doob_heuristic_1949}. This result
 shows
  that the error process $\{\varepsilon_t:t \in [0, 1]\}$ with covariance kernel  \eqref{triangular-kernel} can be represented by
$$\varepsilon_t= \varepsilon(t) = v(t) W(q(t)),$$
where $W = \{W(s): s \in [q(0),q(1) ]\}$ is a Brownian motion on the interval $[q(0), q(1)]$. We use this relationship to represent the process $\{Y_t: t\in [0, 1]\}$ as
\begin{eqnarray*}
Y_t &=& f(t) + \varepsilon_t  =  f(t) + v(t) W(q(t))  , \quad t \in [0, 1] .
\end{eqnarray*}
Dividing by $v(t)$ and using the transformation $s=q(t)$, we get the transformed model
\begin{eqnarray*}
Z_s = g(s) + W(s), \quad s \in [q(0), q(1)],
\end{eqnarray*}
where
$$Z_s= \frac{Y_{q^{-1}(s)}}{v(q^{-1}(s))} \quad \mbox{and} \quad g(s) = \frac{f(q^{-1}(s))}{v(q^{-1}(s))}. $$

Consequently, we obtain for arbitrary $0 \leq t_{i-1}<t_i\leq 1$
\begin{eqnarray*}
\int_{t_{i-1}}^{t_i} \Big(\frac{d}{dt} \left[ \frac{\Phi^{(J)}(t)}{v(t)} \right] \big(\frac{d}{dt}q(t) \big)^{-1} - \mu_i \Big)d\Big(\frac{Y_t}{v(t)} \Big)
&=&  \int_{q(t_{i-1})}^{q(t_i)} \Big(\frac{d}{ds} \left[ \frac{\Phi^{(J)}(q^{-1}(s))}{v(q^{-1}(s))} \right]  - \mu_i  \Big)d\Big(\frac{Y_{q^{-1}(s)}}{v(q^{-1}(s))}\Big) \\
&=& \int_{q(t_{i-1})}^{q(t_i)}\Big( \frac{d}{ds} \tilde\Phi^{(J)}(s) - \mu_i\Big) d Z_s \\
&=& \int_{q(t_{i-1})}^{q(t_i)}\Big( \frac{d}{ds} \tilde\Phi^{(J)}(s) - \mu_i\Big) \big(dg(s) + dW(s)\big),
\end{eqnarray*}
where the function $\tilde\Phi^{(J)}(s)$ is given by $\tilde\Phi^{(J)}(s)= \frac{\Phi^{(J)}(q^{-1}(s))}{v(q^{-1}(s))}$ and we set $t=q^{-1}(s)$. This gives for the transformed derivatives
\begin{eqnarray*}
\frac{d}{dt} \left[ \frac{\Phi^{(J)}(t)}{v(t)} \right] &=& \frac{d}{ds} \left[\frac{\Phi^{(J)}(q^{-1}(s))}{v(q^{-1}(s))}\right] \,  \frac{ds}{dt} = \frac{d}{ds} \left[ \frac{\Phi^{(J)}(q^{-1}(s))}{v(q^{-1}(s))} \right] \,\Big( \frac{d}{ds} q^{-1}(s)\Big)^{-1}, \\
\frac{d}{dt}q(t) &=& \frac{d}{ds} q(q^{-1}(s)) \, \frac{ds}{dt} = \big(\frac{d}{ds} q^{-1}(s) \big)^{-1}.
\end{eqnarray*}
We now introduce the notation
\begin{eqnarray*}
X_i &=& \int_{q(t_{i-1})}^{q(t_i)}\Big( \frac{d}{ds} \tilde\Phi^{(J)}(s) - \mu_i\Big) \big(dg(s) + dW(s)\big)  \quad i=2, \ldots, n.
\end{eqnarray*}
As $W$ is a Brownian motion, the random variables $X_2, \ldots, X_n$ are independent  and  the expected value in \eqref{expec_value} can be rewritten as
\begin{equation}\label{expec_valueXi}
\mathbb{E}\Big[\sum_{i=2}^{n}X_i\sum_{i=2}^{n}X^T_i \Big]= \sum_{i=2}^{n}\mathbb{E}[(X_i- \mathbb{E}[X_i])(X_i- \mathbb{E}[X_i])^T ] + \sum_{i=2}^{n}\mathbb{E}[X_i] \sum_{i=2}^{n}\mathbb{E}[X^T_i].
\end{equation}
Obviously
\begin{eqnarray*}
\mathbb{E}\big[X_i] &=& \int_{q(t_{i-1})}^{q(t_i)} \Big( \frac{d}{ds} \tilde\Phi^{(J)}(s) - \mu_i\Big) \frac{d}{ds}g(s) ds \\
&=&\int_{t_{i-1}}^{t_i} \Big({\frac{d}{dt}\Big[\frac{\Phi^{(J)}(t)}{v(t)}\Big]}  \Big(\frac{d}{dt}{q(t) }\Big)^{-1}  - \mu_i\Big)\Big(\frac{d}{dt} \left[ \frac{{f(t)}}{v(t)} \right]\Big)dt ,
\end{eqnarray*}
and It\^o's isometry gives
\begin{eqnarray*}
&& \mathbb{E}[(X_i- \mathbb{E}[X_i])(X_i- \mathbb{E}[X_i])^T ]\\
&=&\int_{q(t_{i-1})}^{q(t_i)} \Big( \frac{d}{ds} \tilde\Phi^{(J)}(s) - \mu_i\Big)\Big( \frac{d}{ds} \tilde\Phi^{(J)}(s) - \mu_i\Big)^Tds \\
 &=& \int_{t_{i-1}}^{t_i} \Big({\frac{d}{dt}\Big[\frac{\Phi^{(J)}(t)}{v(t)}\Big]}  \Big(\frac{d}{dt}{q(t) }\Big)^{-1}  - \mu_i\Big)\Big({\frac{d}{dt}\Big[\frac{\Phi^{(J)}(t)}{v(t)}\Big]}  \Big(\frac{d}{dt}{q(t) }\Big)^{-1}  - \mu_i\Big)^T \frac{d}{dt}{q(t) }dt.
\end{eqnarray*}
Inserting these representations in \eqref{expec_valueXi} results in \eqref{reformulated_mse}, which proves Proposition \ref{lem_reformulated_mse}.

\end{document}